\documentclass[titlepage,twoside,12pt]{article}
\usepackage{amssymb}
\usepackage{amsfonts}
\begin{document}

{\LARGE \bf Which are the Maximal Ideals ?} \\ \\

{\bf Elem\'{e}r E ~Rosinger} \\ \\
{\small \it Department of Mathematics \\ and Applied Mathematics} \\
{\small \it University of Pretoria} \\
{\small \it Pretoria} \\
{\small \it 0002 South Africa} \\
{\small \it eerosinger@hotmail.com} \\ \\

{\bf Abstract} \\

Ideals of continuous functions which satisfy an {\it off diagonality} condition proved to be
important connected with the solution of large classes of nonlinear PDEs, and more recently,
in General Relativity and Quantum Gravity. {\it Maximal ideals} within those which satisfy
that off diagonality condition are important since they lead to differential algebras of
generalized functions which can handle the largest classes of {\it singularities}. The problem
of finding such maximal ideals satisfying the off diagonality condition is formulated within
some background detail, and commented upon. \\ \\

{\bf 1. The Problem} \\

As presented in section 4, there is a significant interest, both in Mathematics and Physics,
to find out the structure of {\it maximal ideals} ${\cal I}$ in the algebra $( {\cal C} (
\mathbb{R}^n ) )^\Lambda$, among all those ideals which satisfy the {\it off diagonality}
condition \\

(1.1) $~~~ {\cal I} ~\bigcap~ {\cal U}_{\,\Lambda} ( \mathbb{R}^n ) ~=~ \{~ 0 ~\} $ \\

Let us clarify the above notation. First, ${\cal C} ( \mathbb{R}^n )$ is the set of all real
valued continuous functions on $\mathbb{R}^n$, while $\Lambda$ is an arbitrary {\it infinite}
set. Consequently, $( {\cal C} ( \mathbb{R}^n ) )^\Lambda$ is the Cartesian product of
$\Lambda$ copies of ${\cal C} ( \mathbb{R}^n )$, thus it can be identified with ${\cal C}
( \Lambda \times \mathbb{R}^n )$, that is, with the set of real valued continuous functions on
$\Lambda \times \mathbb{R}^n$, where $\Lambda$ is taken with the discrete topology. \\
Clearly, $( {\cal C} ( \mathbb{R}^n ) )^\Lambda$ is a commutative unital algebra over
$\mathbb{R}$, and we have the algebra embedding \\

(1.2) $~~~ {\cal C} ( \mathbb{R}^n ) \ni \psi ~\longmapsto~
                                 u ( \psi ) \in ( {\cal C} ( \mathbb{R}^n ) )^\Lambda $ \\

where $u ( \psi ) = (~ \psi_\lambda ~|~ \lambda \in \Lambda ~)$, with $\psi_\lambda = \psi$,
for $\lambda \in \Lambda$. In this way, the unit element in $( {\cal C} ( \mathbb{R}^n )
)^\Lambda$ is $u ( 1 )$, where $1 \in {\cal C} ( \mathbb{R}^n )$ denotes the constant function
with value 1 defined on $\mathbb{R}^n$. \\
Finally, ${\cal U}_{\,\Lambda} ( \mathbb{R}^n )$ denotes the image of ${\cal C} (
\mathbb{R}^n )$ in $( {\cal C} ( \mathbb{R}^n ) )^\Lambda$ through the algebra embedding (1.2),
thus \\

(1.3) $~~~ {\cal U}_{\,\Lambda} ( \mathbb{R}^n ) ~=~
               \{~ u ( \psi ) ~|~ \psi \in {\cal C} ( \mathbb{R}^n ) ~\} $ \\

is a subalgebra in $( {\cal C} ( \mathbb{R}^n ) )^\Lambda$, and through (1.2), it is
isomorphic with ${\cal C} ( \mathbb{R}^n )$. \\
With the above, the meaning of (1.1) becomes clear, recalling that $\{~ 0 ~\}$ in its right
hand term denotes the trivial zero ideal in $( {\cal C} ( \mathbb{R}^n ) )^\Lambda$. \\

In this way ${\cal U}_{\,\Lambda} ( \mathbb{R}^n )$ is in fact the {\it diagonal} in the
Cartesian product $( {\cal C} ( \mathbb{R}^n ) )^\Lambda$. Thus (1.1) is indeed an {\it off
diagonality} condition on the respective ideals ${\cal I}$ in $( {\cal C} ( \mathbb{R}^n )
)^\Lambda$. \\

As seen in the sequel, the interest in maximal ideals satisfying the off diagonality condition
(1.1) comes from the fact that such ideals lead to the effective construction of differential
algebras of generalized functions which can handle the largest classes of {\it
singularities}. \\

Before going further, let us briefly point to the mathematical nontriviality of the problem in
(1.1). Indeed, as mentioned, $( {\cal C} ( \mathbb{R}^n ) )^\Lambda$ can be identified with
${\cal C} ( \Lambda \times \mathbb{R}^n )$, thus as is well known, Gillman \& Jerison, the
problem of the structure of {\it maximal ideals} ${\cal I}$ in $( {\cal C} ( \mathbb{R}^n )
)^\Lambda$ is closely related to the Stone-\v{C}ech compactification $\beta ( \Lambda \times
\mathbb{R}^n )$ of $\Lambda \times \mathbb{R}^n$, which in itself is a rather involved problem
even in the simplest case of interest above, namely, when $\Lambda = \mathbb{N}$. One of the
reasons which makes $\beta ( \Lambda \times \mathbb{R}^n )$ not easy to deal with is that,
in general, for two completely regular topological spaces $X$ and $Y$, the spaces $\beta ( X
\times Y )$ and $\beta X \times \beta Y$ are different. Furthermore, the space $\beta
\mathbb{N}$ alone is known to be highly nontrivial. \\

On the other hand, in (1.1), one asks the yet more difficult problem of finding the maximal
ideals ${\cal I}$ in $( {\cal C} ( \mathbb{R}^n ) )^\Lambda$ which satisfy the respective {\it
additional condition}, thus they can {\it no longer} be maximal in $( {\cal C} ( \mathbb{R}^n
) )^\Lambda$. Therefore, their structure is quite likely still more complex, Gillman \&
Jerison. \\ \\

{\bf 2. Some Examples of Large Off Diagonal Ideals} \\

As is well known, Gillman \& Jerison, the structure of maximal ideals ${\cal I}$ in the
algebra ${\cal C}( X )$ of real valued continuous functions on a completely regular
topological space $X$ is closely related to certain {\it vanishing} conditions satisfied by
the functions $f \in {\cal I}$. Indeed, in the case $X$ is compact, for instance, then the
maximal ideals ${\cal I}$ in ${\cal C}( X )$ are given by the family of ideals \\

(2.1) $~~~ {\cal M}_p ~=~ \{~ f \in {\cal C} ( X ) ~|~ f ( p ) ~=~ 0 ~\},
                                                        ~~~\mbox{with}~~ p \in X $ \\

In general, when $X$ is not compact, as for instance happens in our case with $X = \Lambda
\times \mathbb{R}^n$, the vanishing condition characterizing the functions in a maximal ideal
is connected with the Stone-\v{C}ech compactification $\beta X$ of $X$. Namely, the maximal
ideals ${\cal I}$ in ${\cal C}( X )$ are given by the family of ideals \\

(2.2) $~~~ {\cal M}^p ~=~ \{~ f \in {\cal C} ( X ) ~|~ p \in cl_{\beta X} Z ( f ) ~\},
                                             ~~~\mbox{with}~~ p \in \beta X $ \\

where $Z ( f ) = \{~ x \in X ~|~ f(x) = 0 ~\}$ is the {\it zero set} of $f$, and $cl_{\beta
X}$ denotes the closure operation in the topology of $\beta X $. \\

In view of the above, we are interested in the case when $X = \Lambda \times \mathbb{R}^n$,
and it is obvious that the ideals ${\cal I}$ in $( {\cal C} ( \mathbb{R}^n ) )^\Lambda$ which
satisfy the off diagonality condition (1.1) must satisfy {\it stronger vanishing} conditions
than those in (2.2), since in view of (1.1), such ideals are significantly smaller than the
maximal ideals in ${\cal C} ( \Lambda \times \mathbb{R}^n )$. \\
Consequently, by looking for such maximal ideals among the ideals satisfying the off
diagonality condition (1.1), we are looking for the {\it weakest} vanishing conditions
satisfied by such ideals. And as seen later, this corresponds to the {\it largest} families of
singularities which the corresponding differential algebras of generalized functions can
handle. \\

A first instance of such stronger vanishing conditions were introduced and used connected with
the so called {\it nowhere dense} ideals, Rosinger [1-19], Mallios \& Rosinger [1-3], Mallios
[1,2], Rosinger \& Walus [1,2], upon which differential algebras of generalized functions were
constructed with the initial aim to solve large classes of nonlinear partial differential
equations. Later, such algebras proved to have a special interest in a variety of basic
theories in Physics, among them General Relativity and Quantum Gravity, Mallios [2]. \\

These nowhere dense ideals are defined as follows. Let $\Lambda = \mathbb{N}$, and let us
denote by ${\cal I}_{nd} ( \mathbb{R}^n )$ the ideal whose elements are all the sequences $w =
( w_0, w_1, w_2, ~.~.~.~ ) \in ( {\cal C} ( \mathbb{R}^n ) )^\mathbb{N}$ of real valued
continuous functions $w_\nu$ which satisfy the {\it asymptotic vanishing} condition \\

(2.3) $~~~ \begin{array}{l}
                   \exists~~~ \Gamma \subset \mathbb{R}^n,~
                                \Gamma ~~\mbox{closed, nowhere dense} ~: \\ \\
                   \forall~~~ x \in \mathbb{R}^n \setminus \Gamma ~: \\ \\
                   \exists~~~ \mu \in \mathbb{N} ~: \\ \\
                   \forall~~~ \nu \in \mathbb{N},~ \nu \geq \mu ~: \\ \\
                   ~~~~~ w_\nu ( x ) ~=~ 0
            \end{array} $ \\

In this case the off diagonality condition (1.1), namely \\

(2.4) $~~~ {\cal I}_{nd} ( \mathbb{R}^n ) ~\bigcap~
                 {\cal U}_{\,\mathbb{N}} ( \mathbb{R}^n ) ~=~ \{~ 0 ~\} $ \\

follows immediately from the fact that the subsets $\mathbb{R}^n \setminus \Gamma$ are always
{\it dense} in $\mathbb{R}^n$, thus a continuous function $\psi \in {\cal C} ( \mathbb{R}^n )$
which vanishes on such a subset must vanish on the whole of $\mathbb{R}^n$. \\
Indeed, let $w = ( w_0, w_1, w_2, ~.~.~.~ ) \in {\cal I} ~\bigcap~ {\cal U}_{\,\mathbb{N}} (
\mathbb{R}^n )$, then in view of (1.3), there exists $\psi \in  {\cal C} ( \mathbb{R}^n )$,
such that $w = u ( \psi )$, therefore $w_\nu = \psi$, for $\nu \in \mathbb{N}$. And then (2.3)
clearly implies that, for a suitable closed and nowhere dense $\Gamma \subset \mathbb{R}^n$,
we have $\psi = 0$ on $\mathbb{R}^n \setminus \Gamma$. Therefore, the continuity of $\psi$
will imply that $\psi = 0$ on the whole of $\mathbb{R}^n$. \\

The meaning of the vanishing condition (2.3) is that the sequences of continuous functions
$w = (~ w_\nu ~|~ \nu \in \mathbb{N} ~)$ in the ideal ${\cal I}_{nd}(\mathbb{R}^n)$ may cover
with their support the {\it singularity} set $\Gamma$, while at the same time, they {\it
vanish asymptotically} outside of it, that is, on $\mathbb{R}^n \setminus \Gamma$. \\
In this way, the ideal ${\cal I}_{nd}(\mathbb{R}^n)$ carries in an {\it algebraic} manner the
information on all the respective sets $\Gamma$ which are sets of {\it singularities} of
generalized functions. And one should recall that such closed and nowhere dense sets $\Gamma$
can have arbitrary {\it large} positive Lebesgue measure, Oxtoby. \\

In view of (2.3), it follows that the sequences of continuous functions in the nowhere dense
ideals ${\cal I}_{nd} ( \mathbb{R}^n )$ satisfy an asymptotic vanishing condition on
corresponding {\it open, dense} subsets of $\mathbb{R}^n$, a condition which is obviously much
stronger than the vanishing conditions in (2.1) or (2.2). \\

As it turned out, however, the nowhere dense ideals ${\cal I}_{nd} ( \mathbb{R}^n )$ were far
from being maximal within the ideals in $( {\cal C} ( \mathbb{R}^n ) )^\mathbb{N}$ which
satisfy the off diagonality condition (1.1). Indeed, in Rosinger [12-16], the following far
larger class of such ideals were introduced and used, see for applications Mallios \& Rosinger
[2,3] and Mallios [2]. \\

Let us consider various families of singularities in $\mathbb{R}^n$, each such family being
given by a corresponding set $\mathcal{S}$ of subsets $\Sigma \subset \mathbb{R}^n$, with each
such subset $\Sigma$ describing a possible set of singularities of a certain given generalized
function. \\

The {\it largest} family of singularities $\Sigma\subset \mathbb{R}^n$ which we can consider
so far is given by \\

(2.5) $~~~ {\cal S}_D ( \mathbb{R}^n ) ~=~
                   \{~ \Sigma\subset \mathbb{R}^n ~~|~~~ \mathbb{R}^n \setminus \Sigma~~~
                          \mbox{is dense in}~~ \mathbb{R}^n ~\} $ \\

And to get an idea how large such singularity sets $\Sigma$ can be, let us note that in the
one dimensional case of $\mathbb{R}$, if we take $\Sigma$ as the set of all {\it irrational}
numbers, then clearly $\Sigma \in {\cal S}_D ( \mathbb{R} )$, since $\mathbb{R} \setminus
\Sigma$ is the set of rational numbers, thus it is dense in $\mathbb{R}$. In this way, it can
happen that a given set $\Sigma$ of singularities has a {\it larger} cardinal than its
complement, that is, than the set of non-singular points. \\

The various families ${\cal S}$ of singularities $\Sigma \subset \mathbb{R}^n$ which we shall
deal with will each satisfy the condition ${\cal S} \subseteq {\cal S}_D ( \mathbb{R}^n )$. \\

Examples of two such families of interest are the following \\

(2.6) $~~~ {\cal S}_{nd}( \mathbb{R}^n ) ~=~ \{~ \Sigma \subset \mathbb{R}^n ~~|~~
                         \Sigma ~~\mbox{closed, nowhere dense in}~ \mathbb{R}^n ~\} $ \\

and \\

(2.7) $~~~ {\cal S}_{Baire~ I} ( \mathbb{R}^n ) ~=~ \{~ \Sigma \subset \mathbb{R}^n ~~|~~
                      \Sigma ~~\mbox{is of first Baire category in}~ \mathbb{R}^n ~\} $ \\

Obviously \\

(2.8) $~~~ {\cal S}_{nd} ( \mathbb{R}^n ) ~\subset~ {\cal S}_{Baire~ I} ( \mathbb{R}^n )
                      ~\subset~ {\cal S}_{D} ( \mathbb{R}^n ) $ \\

And now to the definition of the so called {\it space-time foam} ideals with {\it dense
singularities}, introduced in Rosinger [12-16]. \\

First, let us take any family ${\cal S}$ of singularity sets $\Sigma \subset \mathbb{R}^n$,
family which satisfies the following two conditions \\

(2.9) $~~~ \begin{array}{l}
                  \forall \quad \Sigma \in \mathcal{S}~: \\ \\
                  \quad \quad \mathbb{R}^n \setminus \Sigma
                             \textnormal{ is dense in } \mathbb{R}^n
           \end{array} $ \\

and \\

(2.10) $~~~ \begin{array}{l}
                   \forall \quad \Sigma,~ \Sigma ' \in \mathcal{S}~: \\ \\
                   \exists \quad \Sigma '' \in \mathcal{S}~: \\ \\
                   \quad \quad \Sigma \cup \Sigma ' \subseteq \Sigma ''
            \end{array} $ \\

Clearly, we shall have the inclusion $\mathcal{S} \subseteq \mathcal{S}_{\mathcal{D}} (
\mathbb{R}^n )$ for any such family $\mathcal{S}$. \\

It is easy to see that both families $\mathcal{S}_{nd} (\mathbb{R}^n)$ and
$\mathcal{S}_{Baire~ I} ( \mathbb{R}^n )$ satisfy conditions the (2.9) and (2.10). \\
On the other hand, the family ${\cal S}_{D} ( \mathbb{R}^n )$ as a whole does {\it not}
satisfy condition (2.10). Indeed, one can partition $\mathbb{R}^n$ into two subsets $\Sigma$
and $\Sigma '$, both of which are dense in $\mathbb{R}^n$, thus both belong to ${\cal S}_{D}
( \mathbb{R}^n )$. In this case, in (2.10), we would have to have $\Sigma '' = \mathbb{R}^n$,
which obviously does not belong to ${\cal S}_{D} ( \mathbb{R}^n )$. \\

Now, as the second ingredient, and so far independently of any $\mathcal{S}$ above, we take
any right directed partial order $L=(\Lambda,\leq)$. In other words, $L$ is such that for each
$\lambda,\, \lambda ' \in \Lambda$ there exists $\lambda ''\in \Lambda$ with $\lambda,
\lambda ' \leq \lambda ''$. Here we note that the choice of $L$ may at first appear to be
completely independent of $\mathcal{S}$, yet in certain specific instances the two may be
somewhat related, with the effect that $\Lambda$ may have to be large, see Rosinger [15]. \\

Although we shall only be interested in singularity sets $\Sigma \in \mathcal{S}_{\mathcal{D}}
(\mathbb{R}^n)$, the following {\it ideal} can in fact be defined for any $\Sigma \subseteq
\mathbb{R}^n$. Indeed, let us denote by \\

(2.11) $~~~ {\cal J}_{L,~\Sigma}(\mathbb{R}^n) $ \\

the {\it ideal} in $( {\cal C} (\mathbb{R}^n))^{\Lambda}$ of all the sequences of continuous
functions indexed by $\lambda \in \Lambda$, namely, $w = (~ w_{\lambda} ~|~ \lambda \in
\Lambda ~) \in ( {\cal C} (\mathbb{R}^n))^{\Lambda}$, sequences which {\it outside} of the
singularity set $\Sigma$ will satisfy the {\it asymptotic vanishing} condition \\

(2.12) $~~~ \begin{array}{l}
                 \forall \quad x \in \mathbb{R}^n \setminus \Sigma~: \\ \\
                 \exists \quad \lambda \in \Lambda~: \\ \\
                 \forall \quad \mu \in \Lambda,~ \mu \geq \lambda~: \\ \\
                 \quad \quad w_{\mu}(x)=0
             \end{array} $ \\

This means that the sequences of continuous functions $w = (~ w_\lambda ~|~ \lambda \in
\Lambda ~)$ in the ideal $\mathcal{J}_{L,~\Sigma}(\mathbb{R}^n)$ may {\it cover} with their
support the singularity set $\Sigma$, and at the same time, they {\it vanish asymptotically}
outside of it, that is, on $\mathbb{R}^n \setminus \Gamma$. \\
In this way, the ideal ${\cal J}_{L,~\Sigma} ( \mathbb{R}^n )$ carries in an {\it algebraic}
manner the information on the singularity set $\Sigma$. \\

Here however, with the ideals ${\cal J}_{L,~\Sigma}(\mathbb{R}^n)$, the asymptotic vanishing
happens on $\mathbb{R}^n \setminus \Sigma$. And as seen below, we shall only require that such
sets $\mathbb{R}^n \setminus \Sigma$ be {\it dense} in $\mathbb{R}^n$, in other words that
$\Sigma \in {\cal S}_D ( \mathbb{R}^n )$. Thus the corresponding vanishing conditions are
significantly {\it weaker} than that required in the case of the nowhere dense ideals ${\cal
I}_{nd} ( \mathbb{R}^n )$. \\

In follows that the nowhere dense ideals only allow singularities on closed and nowhere dense
subsets $\Gamma \subset \mathbb{R}^n$, whose complementaries $\mathbb{R}^n \setminus \Gamma$
are therefore open and dense in $\mathbb{R}^n$. On the other hand, the case of the ideals
${\cal J}_{L,~\Sigma} ( \mathbb{R}^n )$ in which we shall be interested, namely when $\Sigma
\in {\cal S}_D ( \mathbb{R}^n )$, allow arbitrary and much {\it larger singularities} $\Sigma
\subset \mathbb{R}^n$, as long their complementaries $\mathbb{R}^n \setminus \Sigma$ are still
dense in $\mathbb{R}^n$. \\

We note that the assumption about $L=(\Lambda,\leq)$ being right directed is used in proving
that ${\cal J}_{L,~\Sigma} ( \mathbb{R}^n )$ is indeed an ideal, more precisely that, for $w,
w' \in {\cal J}_{L,~\Sigma} ( \mathbb{R}^n )$, we have $w+w' \in {\cal J}_{L,~\Sigma} (
\mathbb{R}^n )$. \\

Now, it is easy to see that for $\Sigma,~ \Sigma\, ' \subseteq \mathbb{R}^n$ we have \\

(2.13) $~~~ \Sigma \subseteq \Sigma\, ' ~~\Longrightarrow~~
                         {\cal J}_{L,~\Sigma} ( \mathbb{R}^n ) \subseteq
                                    {\cal J}_{L,~\Sigma\, '} ( \mathbb{R}^n ) $ \\

In this way, for any family ${\cal S}$ of singularity sets $\Sigma \subset \mathbb{R}^n$
satisfying (2.9), (2.10), it follows that \\

(2.14) $~~~ {\cal J}_{L,~{\cal S}} ( \mathbb{R}^n ) ~=~ \bigcup_{\, \Sigma \in {\cal S}}~
                        {\cal J}_{L,~\Sigma} ( \mathbb{R}^n ) $ \\

is also an {\it ideal} in $( {\cal C} ( \mathbb{R}^n))^{\Lambda}$. \\

It is important to note that for suitable choices of the right directed partial orders $L$,
the ideals ${\cal J}_{L,\Sigma} ( \mathbb{R}^n )$, with $\Sigma \in {\cal S}_D ( \mathbb{R}^n
)$, are {\it nontrivial}, that is, they do not reduce to the zero ideal $\{~0~\}$, Rosinger
[15, section 2]. Thus in view of (2.14), the same will hold for the ideals ${\cal
J}_{L,~{\cal S}} ( \mathbb{R}^n )$. \\

Let us conclude by showing that the ideals ${\cal J}_{L,~{\cal S}} ( \mathbb{R}^n )$ satisfy
the off diagonality condition (1.1), namely \\

(2.15) $~~~ {\cal J}_{L,~{\cal S}} ( \mathbb{R}^n ) \bigcap
                 {\cal U}_{\,\Lambda} ( \mathbb{R}^n ) ~=~ \{~ 0 ~\} $ \\

for every family ${\cal S}$ of singularities which satisfies (2.9) and (2.10). Indeed, let \\

$~~~ w = (~ w_{\lambda} ~|~ \lambda \in \Lambda ~) \in
                     {\cal J}_{L,~{\cal S}} ( \mathbb{R}^n ) \bigcap
                                       {\cal U}_{\,\Lambda} ( \mathbb{R}^n ) $ \\

then in view of (2.14), there exists $\Sigma \in {\cal S}$, such that \\

$~~~ w = (~ w_{\lambda} ~|~ \lambda \in \Lambda ~) \in
                          {\cal J}_{L,~\Sigma} ( \mathbb{R}^n ) $ \\

On the other hand, we have $w = u ( \psi )$, for a certain $\psi \in {\cal C} ( \mathbb{R}^n
)$. Thus $w_\lambda = \psi$, for $\lambda \in \Lambda$. And then (2.12) implies that $\psi =
0$ on $\mathbb{R}^n \setminus \Sigma$, which means that $\psi = 0$ on the whole of
$\mathbb{R}^n$, since $\mathbb{R}^n \setminus \Sigma$ is dense in $\mathbb{R}^n$, in view of
(2.9). \\ \\

{\bf 3. Towards Finding the Maximal Ideals} \\

Let $X$ be a completely regular topological space. The relations (2.1) and (2.2) indicate that
important properties of the algebra ${\cal C} ( X )$ of real valued continuous functions on
$X$ may be expressed in terms of the topology of $X$ or of $\beta X$. Such a translation of
algebraic properties into topological ones, or vice versa, may often prove useful, and we
shall mention some of them. \\

{\bf Customary Notations and Facts.} We recall, Gillman \& Jerison, the following notations.
First \\

(3.1) $~~~ {\bf Z} ( X ) ~=~ \{~ Z ( f ) ~~|~~ f \in {\cal C} ( X ) ~\} $ \\

is the set of all {\it zero-sets} of functions in ${\cal C} ( X )$. Each such $Z ( f )$ is
closed in $X$, since the respective $f$ are continuous. However, ${\bf Z} ( X )$ need not in
general be the set of all closed subsets of $X$. \\

A family ${\cal F}$ of zero-sets is called a {\it z-filter} on $X$, if and only if it
satisfies the following conditions \\

(3.2) $~~~ \begin{array}{l}
               ~~~~~*) \quad \phi \notin {\cal F} \\ \\
               ~~**) \quad Z,~ Z\,' \in {\cal F} ~~\Longrightarrow~~
                                             Z \bigcap Z\,' \in {\cal F} \\ \\
               ~***) \quad Z \in {\cal F},~ Z\,' \in Z ( X ),~ Z \subseteq Z\,'
                                       ~~\Longrightarrow~~ Z\,' \in {\cal F}
           \end{array} $ \\ \\

For any ideal ${\cal I}$ in the algebra ${\cal C} ( X )$ of real valued continuous functions
on $X$, we denote \\

(3.3) $~~~ {\bf Z} ( {\cal I} ) ~=~ \{~ Z ( f ) ~~|~~ f \in {\cal I} ~\} $ \\

This is known to be a {\it z-filter} on $X$. \\

Further, for every z-filter ${\cal F}$ on $X$, we denote \\

(3.4) $~~~ {\bf I} ( {\cal F} ) ~=~ \{~ f \in {\cal C} ( X ) ~~|~~
                                            Z ( f ) \in {\cal F} ~\} $ \\

This is known to be an ideal in ${\cal C} ( X )$. Furthermore \\

(3.5) $~~~ {\bf Z} ( {\bf I} ( {\cal F} ) ) ~=~ {\cal F},~~~
           {\bf I} ( {\bf Z} ( {\cal I} ) ) ~\supseteq~ {\cal I} $ \\

{\bf Sets and Ideals.} With the above customary notations, and based on (2.2), for any ideal
${\cal I}$ in the algebra ${\cal C} ( X )$ of real valued continuous functions on $X$, we
denote \\

(3.6) $~~~ P ( {\cal I} ) ~=~ \{~ p \in \beta X ~~|~~
                          {\cal I} ~\subseteq~ {\cal M}^p ~\} ~\subseteq~ \beta X $ \\

Clearly $P ( {\cal I} ) \neq \phi$, since every ideal ${\cal I}$ in ${\cal C} ( X )$ is
contained in some maximal ideal ${\cal M}^p$ in ${\cal C} ( X )$. Also obviously \\

(3.7) $~~~ {\cal I} ~\subseteq~ \bigcap_{p \in P({\cal I})} {\cal M}^p $ \\

however, equality need not always hold. Indeed, let $X = [ -1, 1 ] \subset \mathbb{R}$, and
${\cal I}$ be the principal ideal in ${\cal C} ( X )$ generated by the function $id_X$, that
is \\

$~~~ {\cal I} = \{ ~ f \in {\cal C} ( X ) ~~|~~
            \exists~ g \in {\cal C} ( X ) ~:~ f ( x ) =
                                           x g ( x ),~~ x \in X ~\} $ \\

In this case $P ( {\cal I} ) = \{ 0 \}$, while ${\cal I} \neq {\cal M}^0$, since if we define
$f \in {\cal C} ( X )$ by $f ( x ) = \sqrt |~ x ~|$, for $x \in X$, then $f \in {\cal M}^0$,
but $f \notin {\cal I}$. \\

Further, (3.6) and (2.2) give \\

(3.8) $~~~ P ( {\cal I} ) ~\subseteq~ \bigcap_{f \in {\cal I}}~ cl_{\beta X} Z ( f ) $ \\

Based on the above, and as an extension of (2.2), for every subset $A \subseteq \beta X$, let
us define the ideal in ${\cal C}( X )$ given by \\

(3.9) $~~~ {\cal I}^A ~=~ \bigcap_{p \in A} {\cal M}^p $ \\

Then obviously ${\cal I}^A = {\cal M}^p$, for $A = \{~ p ~\}$, with $p \in \beta X$. Further,
we have \\

(3.10) $~~~ A ~\subseteq~ P ( {\cal I}^A )~~~ \mbox{for}~~ A ~\subseteq~ \beta X $ \\

Also \\

(3.11) $~~~ {\cal I}^A ~=~ \{ ~ f \in {\cal C} ( X ) ~~|~~
                        A ~\subseteq~ cl_{\beta X} Z ( f ) ~\} $ \\

and \\

(3.12) $~~~ {\cal I} ~\subseteq~ {\cal I}^{P ( {\cal I} )} $ \\

{\bf The Problem in Section 1.} In that case we have $X = \Lambda \times \mathbb{R}^n$, thus
we obtain \\

(3.13) $~~~ Z ( u ( \psi ) ) ~=~ \Lambda \times Z ( \psi )~~~
                                      \mbox{for}~~ \psi \in {\cal C} ( X ) $ \\

where $Z$ in the left hand term is considered in $X = \Lambda \times \mathbb{R}^n$, while $Z$
in the right hand term is defined in $\mathbb{R}^n$. \\

Now for convenience, let us first consider the following {\it particular} case of the problem
in section 1, namely : \\

Which are the subsets $A \subseteq \beta X$, such that \\

(3.14) $~~~ {\cal I}^A ~\bigcap~ {\cal U}_{\,\Lambda} ( \mathbb{R}^n ) ~=~ \{~ 0 ~\}~~~ ? $ \\

In other words, which are the subsets $A \subseteq \beta X$, such that, given any $\psi \in
{\cal C} ( \mathbb{R}^n )$, we have \\

(3.15) $~~~ u ( \psi ) \in {\cal I}^A ~~\Longrightarrow~~
                                     \psi ~=~ 0 ~~~\mbox{on}~~ \mathbb{R}^n~~~ ? $ \\

In view of (3.11), (3.13), the problem (3.15) is equivalent with finding the subsets $A
\subseteq \beta X$, such that, given any $\psi \in {\cal C} ( \mathbb{R}^n )$, we have \\

(3.16) $~~~ A ~\subseteq~ cl_{\beta ( \Lambda \times \mathbb{R}^n )}
              ( \Lambda \times Z ( \psi ) )
                  ~~\Longrightarrow~~ \psi ~=~ 0~~~ \mbox{on}~~ \mathbb{R}^n $ \\

Clearly, the larger $A$, the more likely that $Z ( \psi )$ is large, thus the implication in
(3.16) may hold. An obvious {\it sufficient} condition for (3.16) is the following \\

(3.17) $~~~ pr_{_{\mathbb{R}^n}}~ ( A ~\bigcap~ ( \Lambda \times \mathbb{R}^n ) )~~~
                                                     \mbox{dense in}~~ \mathbb{R}^n $ \\

{\bf The Case of Space-Time Foam Ideals.} Let us note here what happens in the case of
sequences of continuous functions, see (2.11),
(2.12) \\

(3.18) $~~~ w = (~ w_{\lambda} ~|~ \lambda \in \Lambda ~) \in
                                        {\cal J}_{L,~\Sigma}(\mathbb{R}^n) $ \\

where $\Sigma \subset \mathbb{R}^n$ is such that $\mathbb{R}^n \setminus \Sigma$ is dense in
$\mathbb{R}^n$. In view of (2.12), we have \\

(3.19) $~~~ Z ( w ) ~=~ \bigcup_{x \in \mathbb{R}^n \setminus \Sigma}~ (~ [ \lambda_x >
                     \times \{~ x ~\} ~) ~\subseteq~ \Lambda \times \mathbb{R}^n ~=~ X $ \\

where $\lambda_x \in \Lambda$ and $[ \lambda_x >\, = \,\{~ \mu \in \Lambda ~|~ \mu \geq
\lambda_x ~\}$. \\

Obviously, for the set $A = Z ( w )$ in (3.19), we have \\

(3.20) $~~~ A~~~ \mbox{satisfies}~~ (3.17) ~~~\Longleftrightarrow~~~
            \mathbb{R}^n \setminus \Sigma~~~ \mbox{dense in}~~ \mathbb{R}^n $ \\

which clarifies the conditions in (2.5) and (2.9). \\ \\

{\bf 4. On the Off Diagonality Condition} \\

In Rosinger [6, chap. 3, pp. 65-119] was for the first time given an {\it algebraic
characterization} of all differential algebras $A$ of generalized functions which contain the
Schwartz distributions ${\cal D}'$. Based on this characterization, the {\it infinite} class
of all differential algebras of generalized functions containing the Schwartz distributions
was constructed in Rosinger [6], see also Rosinger [1-5,7-19]. \\

These differential algebras prove to be of particular interest in the solution of large
classes of nonlinear PDEs, see Rosinger [1-19], Rosinger \& Walus [1,2], Colombeau, Biagioni,
Oberguggenberger, Grosser et.al., and the literature cited there, as well as the subject 46F30
in the AMS Subject Classification, at www.ams.org/msc/46Fxx.html \\

Recently, a class of these differential algebras, called {\it space-time foam algebras}\, was
introduced in Rosinger [12-16]. This class proves to be of special interest in setting up a
Differential Geometry which can handle the largest classes of {\it singularities} so far in
the literature. This property is of special interest in General Relativity and Quantum Gravity,
see Mallios \& Rosinger [1-3] and Mallios [1,2]. \\

Returning to the mentioned algebraic characterization, it also shows that the Colombeau
algebras, see Colombeau, Biagioni, Oberguggenberger, are a {\it particular} case of the
infinite class of differential algebras constructed in Rosinger [1-19], see also the comment
in Grosser et.al. [p. 7]. \\

In some more detail the situation is as follows. The differential algebras $A$ of generalized
functions in the mentioned infinite class are of the {\it quotient} form \\

(4.1) $~~~ A ~=~ {\cal A} / {\cal I} $ \\

where \\

(4.2) $~~~ {\cal I} ~\subset~ {\cal A} ~\subseteq~
                       ( {\cal C}^\infty ( \mathbb{R}^n ) )^\Lambda $ \\

with $\Lambda$ being suitable infinite sets, while ${\cal A}$ are subalgebras in $( {\cal
C}^\infty ( \mathbb{R}^n ) )^\Lambda$, and ${\cal I}$ are ideals in ${\cal A}$. \\

In this case, the {\it off diagonality} condition which characterizes such quotient algebras
$A$ when ${\cal A} = ( {\cal C}^\infty ( \mathbb{R}^n ) )^\Lambda$, is given by \\

(4.3) $~~~ {\cal I} ~\bigcap~ {\cal U}\,^\infty_\Lambda ( \mathbb{R}^n ) ~=~ \{~ 0 ~\} $ \\

where ${\cal U}\,^\infty_\Lambda ( \mathbb{R}^n )$ is the diagonal in the Cartesian product
$( {\cal C}^\infty ( \mathbb{R}^n ) )^\Lambda$. \\

However, at a first approach, one may set aside the ${\cal C}^\infty$-smoothness involved in
(4.1) - (4.3), and instead, consider the more general and merely continuous setup presented in
section 1. \\


\begin{thebibliography}{99}
\bibitem{} Rosinger E E [17] : Can there be a general nonlinear PDE theory for the existence
of solutions ? math.AP/0407026
\end{thebibliography}
\end{document}